\documentclass[12pt]{article}
\usepackage{amssymb,amsmath,amsthm}

 \usepackage{latexsym}

\title{\bf \boldmath }

\newtheorem{theorem}{Theorem}[section]
\newtheorem*{thma}{Theorem}

\newtheorem{definition}[theorem]{Definition}
\newtheorem{lemma}[theorem]{Lemma}

\newtheorem{corollary}[theorem]{Corollary}

\newcommand{\Fitt}{{\rm Fitt}}
\newcommand{\po}{{\phi_\pi}}

\newcommand{\N}{{\mathbb{N}}}
\newcommand{\Z}{{\mathbb{Z}}}
\newcommand{\Ker}{{\rm Ker}}
\newcommand{\upi}{{\cal U}_\pi}

\newcommand{\mapright}[1]{{\smash{\mathop{\longrightarrow}\limits^{#1}}}}

\newcommand{\rmapdown}[1]{{\Big\downarrow\rlap{$\vcenter{\hbox{$\scriptstyle#1$}}$}}}
\newcommand{\lmapdown}[1]{{\vcenter{\hbox{$\scriptstyle#1$}}\Big\downarrow
               \vcenter{\hbox{$\phantom{\scriptstyle#1}$}}}}
\newcommand{\mapdown}{{\Big\downarrow}}

\newcommand{\z}{{\Z_\pi}}

\newcommand{\comdiag}[1]{\begin{matrix}#1\end{matrix}}
\newcommand{\ppi}{{\cal P}_\pi}
\newenvironment{pf}{\noindent{\bf Proof.}}{\quad\medskip}

\begin{document}
\bibliographystyle{plain}
\title{\bf \boldmath Localizations of finitely generated soluble groups}
\author{ Niamh O'Sullivan\\
School of Mathematical Sciences, \\Dublin City University,\\  Dublin  9, Ireland \\ email:niamh.osullivan@dcu.ie}
\date{}

\maketitle

\begin{abstract}
The study of localizations of groups has concentrated on group
theoretic properties which are preserved by localization. In this
paper we look at finitely generated soluble groups and determine
when the local groups associated with them are soluble.
\end{abstract}

\section{Introduction}  Let $G$ be a
group and let $\pi$ be a set of primes. The group $G$ is said to
be $\pi$-local, if the map $x \mapsto x^q $ is bijective, for
every prime $q\not \in \pi $. A $\pi$-localization of any group
$G$ is a homomorphism $\po :G \rightarrow G_\pi, $ where $G_\pi$
is a $\pi$-local group, with the universal property that given any
homomorphism $\theta:G \rightarrow H,$ where $H$ is $\pi-$local,
there exists a unique homomorphism $\theta_\pi :G_{\pi}
\rightarrow H$ such that $\theta_\pi \po =\theta .$ Localizations
always exist and are unique up to isomorphism \cite{rib}. Throughout this paper we shall assume that $\pi$ is not the set of all primes. Given any set of primes $\pi$, $\pi'$
 will denote its complement in the set of all primes. If $m$ is an
 integer whose prime divisors all lie in $\pi'$, then we write
 $m\in\pi' .$ We say that $x$ is a $q$-th root of  $y,$
if $x^q=y.$ A group $G$ is a $\upi$-group   if the map $x \mapsto x^m  $ is injective whenever $m\in \pi'.$

 \par Localizations have been studied extensively, with particular focus
on properties which are preserved by localization. For example,
Ribenboim \cite{rib} and Casacuberta and Castellet
\cite{casacast}, have shown that localization preserves nilpotency
and  virtual nilpotency respectively. It is therefore natural to
ask whether localization preserves solubility. In general it does
not: we have shown in
 \cite{nos4}, that if $G$ is a free soluble  group of derived length greater than $1$, then $G_\pi$ is not soluble. The aim of
 this paper is to tackle the question for finitely generated  soluble groups.

\par The study of localizations of soluble groups is intrinsically
more
 difficult than that of nilpotent groups.   Localization preserves surjectivity, Proposition 10 \cite{rib}, but  does not preserve injectivity. Localization is an exact functor on the class of
nilpotent groups, Proposition 18 \cite{rib} but is not an exact functor on the
class of soluble groups.  For example, when $\pi =\{3\} ,$ we have
the following diagram:

$$\comdiag { 1 &\mapright{}& C_3&\mapright{}& S_3 &\mapright{}& C_2&\mapright{}& 1 \cr &&\lmapdown{\po }
&&\rmapdown{\po} &&\rmapdown{\po} &&\cr  && C_3 && 1  &\mapright{}&
1&\mapright{}& 1 .}$$
 Due to this, very little is known about the structure of $G_\pi$
 when $G$ is a soluble group.

 \par   Baumslag \cite{baum1}
  defined a class of groups  which was later modified by Cassidy \cite{cass}. As Cassidy's class of groups does not include
   groups  with elements of  order $2,$ we shall make a slight modification
   and this new class of groups will be called   ${\cal P}_\pi$. If $G\in \ppi, $ then,  using Baumslag's methods, we can
   construct $G_\pi$ explicitly. We prove the following
\begin{thma}   Let  $G \in \ppi$ be   a  soluble group which is not $\pi$-local.
Then $G_\pi $ is not soluble. \end{thma}
\par As localization preserves
surjectivity, an immediate consequence of this is that  if  $G$ has a quotient which lies  in  $\ppi$, then $G_\pi$ will not be soluble.
We determine when a finitely generated soluble group has such a quotient and prove the following.

\begin{thma} Let $G$ be a finitely generated metabelian group and
let $\pi$ be a finite set of primes. Then
$G$ either has a normal subgroup $N$ such that $G/N\in \ppi$ and so $G_\pi$ is  not soluble or  $G_\pi$ is
virtually nilpotent.\end{thma}
\begin{thma}
Let  $G$ be a finitely generated soluble group.
Then either $G_\pi$ is virtually nilpotent, for all sets of primes
$\pi ,$ or there exists a finite set of primes $\tau$ such that
$G_\pi$ is not soluble whenever $\tau\subseteq \pi .$
\end{thma}

\section {The class ${\cal P}_\pi.$}

  Baumslag\cite{baum1} and  Cassidy \cite{cass}  each  defined  a class of
   groups, $\ppi$, with the property  that if $G\in \ppi,$ then $G$   embeds in $G_\pi$.  If $G\in \ppi$ and $g\in G,$ then
   the idea
   is to embed $G$ in a group $K\in \ppi$ such that  $g$ has $p$-th roots in $K$, for all $p\in \pi' .$ This  is done
   repeatedly until $G$ embeds in  a group $\overline G \in \ppi$ having $p$-th roots for all the elements of $G$
    and for all $p\in \pi' .$
One then repeats the process on the group $\overline G$ and on the
resulting groups and hence obtains an ascending sequence of groups
in  $\ppi .$  The union of this ascending sequence is   a
$\pi$-local group containing $G$ and in fact is $G_\pi.$  We shall
make slight modifications to the definition of $\ppi$ and show
that we can also construct $G_\pi$ in a similar manner if $G\in
\ppi.$ We will show in Theorem 2.6 that if $G\in \ppi$ is a
soluble group, then the group, that we've called $K$ above,
obtained at the first step is not soluble and hence $G_\pi$ is not
soluble.

\begin{definition}  Let $\pi$ be a set of primes. A group $G$ is said to
be in the class ${\cal P}_\pi$ if  \begin{itemize} \item[(i)] $G\in\upi$;

 \item[(ii)]  if   $ x\in G$ has no $m$-th root, for some $m\in \pi',$ and $$D_G (x)= \{ g\in
G\mid C_G (x)\cap C_G (x)^g\not=1\},$$ then  $C_G (x)$ is
torsion-free and   one of the following holds:
\begin{itemize} \item[(a)] $D_G (x)=C_G (x)$ so that  $C_G (x)$ is  malnormal;
\item[(b)] $|D_G (x):C_G(x)|=2$ and $C_G(x) \le \z .$
\end{itemize}\end{itemize}
\end{definition}
 \vskip .5cm \par If $G \in \ppi $ and $ x\in G$ has no $m$-th
 root, for some $m\in \pi',$ then we adjoin $p$-th roots of $x,$
 for all $p\in \pi',$ by making slight modifications to the method used in \cite{cass}. Let $D=D_G (x)$ and let $C=C_G (x).$ If
 $D=C,$  then set $R=\{ g\in C\mid  g=x^{a/b} ,\; b\in \pi'\}$,
 this is a central subgroup of $C$ and is isomorphic to a subgroup
 of $\z$ through the monomorphism $i:R\rightarrow \z$ that
sends
  $x$ to $ 1 .$ Let $P$ denote the central product of $C$ and  $\z $ over
$R.$  So that $P$ is the quotient $(C\times \z)/A$ where $A=\{ (y
,i (y^{-1})) \mid y\in R \} .$ If $|D:C|=2 ,$ then, without loss of generality, we can assume
that we've chosen $x$ so that we can view $C$ as a
subgroup of $\z$ through the monomorphism $i:C\rightarrow \z$ that
sends
  $x$ to $ 1 .$ The action of $D$ on $C$ induces an action of $D$
  on $\z$ and we let  $P$ denote the semi-direct product of $D$ and  $\z $ over
$C.$  So that   $P$ is the quotient $(D\ltimes \z)/A$ where $A=\{
(y ,i (y^{-1})) \mid y\in C \} .$ In both cases let $K$ denote the
generalized free product of $G$ and $P$ with $D$ amalgamated, that
is
$$K= G*_D P
$$ (without loss of generality we shall now refer to $D\le  P$).
We can choose  a left transversal, $S_G,$ to $D$ in $G$ and a left
transversal, $S_P,$ to $D$ in $P.$ Then every element $g\in K$ can
be written uniquely in normal form: $$g=s_1s_2\dots s_nz$$ where
each $s_i$ belongs to either $S_G$ or $S_P,$   successive
components $s_i$ and $s_{i+1}$ do not belong to the same
transversal $S_G$ or $S_P$ and the last component $z$ belongs to
the common subgroup $D.$ We call $n$ the length of the normal form
and it is usually denoted $\lambda (g).$  In order to prove that $K\in
\ppi ,$ we only need to adapt the work of Baumslag \cite{baum1}
and Cassidy \cite{cass} slightly to take
 into account condition 2.2 (ii) b so that we will mainly
 reference their work here rather than replicate it.

 \begin{lemma} Let $P$ be either of the groups described above.
 Then $P\in\upi$   and if $|D:C|=2,$ then
 $P$ is a $\pi$-local group. \end{lemma}
 \begin{pf} If $D=C,$ then, by Lemma 3 \cite{cass}, $P\in\upi .$  If $|D:C|=2$ and $C \le \z$,  then $\z\triangleleft P$
 and $P/\z\cong C_2 .$ As $G\in\upi ,$  $2\in \pi
 .$ Therefore $P$ is an extension of a $\pi$-local group by a
 finite $\pi$-group and hence is $\pi$-local, by Theorem 11.5 \cite{baum1}. $\square$
 \end{pf}

\begin{lemma} Let $a\in K.$
\begin{itemize} \item[(i)] If $a\in P$ and $a^2\not= 1,$ then $ C_K (a)\le P.$
 Furthermore if $a$ does not have all its $\pi'$-roots, then
$C_K (a)=P.$
 \item[(ii)] If $a\in G$ and $a$ is not conjugate to an element of
$C$, then $C_G (a)=C_K (a).$  \item[(iii)] If $a$ is conjugate to
neither an element of $G$ nor an element of $P$, then $C_K (a)$ is
an infinite cyclic group. Moreover if $a$ is cyclically reduced,
then every non-trivial $y\in C_K (a)$ is cyclically reduced and is
conjugate to neither an element of $G$ nor an element of $P$.

\end{itemize}
\end{lemma}
\begin{pf} This follows directly from Lemmas 5-8 \cite{cass}.
We note that Lemma 8 only applies to the case when $D=C$ as when
$|D:C|=2,$ $P$ is $\pi$-local by Lemma 2.2. $\square$

\end{pf}

\begin{theorem} Let $\displaystyle K= G*_D P$. Then   $K\in \ppi$.
\end{theorem}
\begin{pf} The fact that $K\in \upi $  follows immediately from Lemmas 2.2 and 2.3 as in
Theorem 29.1 p. 266 \cite{baum1}. Suppose that    $a\in K$ and $a$ does
not have a $p$-th root in $K,$ for some $p\in \pi'$.  Since  $K\in\upi$, $a$ has infinite order. If $a\in P$,  then $P$ is not
$\pi$-local, $C=D$ and, by Lemma 2.3 (i), $C_K (a)= P$ which is torsion-free.  In the other cases we see, by Lemma 2.3 (ii) and (iii), that
  $C_K (a)$ is torsion-free.
\par$\;$\par    We must also show that either $D_K (a)=C_K (a)$ or $|D_K ( a):C_K(a)|=2$ and $C_K (a) \le \z.$  As in Proposition 2 \cite{cass}, we can assume
that $a$ is cyclically reduced. Suppose first that $a \in P. $  Then we've seen above that $C=D$ and $C_K (a)=P$ and  as in Proposition 2
\cite{cass}, $C_K (a)$ is malnormal so $D_K (a)=C_K (a).$
 Now suppose that $a \in G$ but is not conjugate to an element
of $D.$  Then $C_K (a)= C_G (a)$, by Lemma 2.3 (ii). As in
Proposition 2 \cite{cass}, $C_K (a) \cap C=1$ and since $C_G (a)$
is torsion-free this implies that $C_K (a)\cap D=1.$  Following
the proof of Proposition 2 \cite{cass} we see that $D_G (a)=D_K
(a).$ Finally suppose that $\lambda (a)\ge 2.$ Then $C_K (a)$ is
an infinite cyclic group and, by Proposition 2 \cite{cass}, $x\in
D_K (a)$ implies that either $x\in C_K (a)$ or $x^2=1$ and
$a^x=a^{-1}.$  So that either $D_K (a)=C_K (a)$ or $|D_K (a):C_K(a)|=2$ and $C_K (a) \le \z.$ Therefore $K\in \ppi.$ $\square$

\end{pf}
\begin{theorem} Let $G\in \ppi$. Then $G$ can be embedded in
$G_\pi$
\end{theorem}
\begin{pf} We use an
classical ascending tower argument to embed $G$ in $G_\pi$  as
stated at the beginning of this section and as in  Theorem 1
\cite{cass}. $\square$ \end{pf}

\begin{theorem}   Let  $G \in \ppi$ be   a  soluble group which is not $\pi$-local.
Then $G_\pi $ is not soluble. \end{theorem}
\begin{pf}  If $x\in \Fitt (G),$ the
Fitting subgroup of $G$, then $x\in N$ where $N$ is a nilpotent
normal subgroup of $G$.   Therefore $Z(N) \le C_G (x)$ and $D_G
(x)= G.$
  As $x$ does not satisfy 2.1 (ii), $\pi'$-roots
for $x$ exist in $G.$  If $G$ is virtually nilpotent, then $\Fitt (G)$ is of finite $\pi$-index
in $G$  and hence  $G$ is $\pi$-local, by Corollary 15.3  and Theorem 11.5 \cite{baum1}.
\par $\;$
\par We can therefore  assume that  $G$ is not virtually nilpotent and that there exists $x\in
G\backslash \Fitt (G),$ $m\in \pi'$ such that $x$ has no $m$-th
root, $C=C_G (x)$ is torsion-free and $D=D_G (x)$ satisfies 2.2
(ii). We adjoin $\pi'$-roots for $x$ by setting
$$K= G*_C P $$ as before.
Then, by  Theorems 2.4 and 2.5, $K$ embeds in $K_\pi \cong G_\pi
.$ By Proposition 11.22 \cite{ls}, $K$ has a  free subgroup of rank $2$ and hence $K$ is not soluble. Therefore $G_\pi$
is not soluble. $\square$\end{pf}

\section{Soluble groups}  We now look at finitely generated soluble groups in
general. We reiterate  that we are assuming throughout that  $\pi$ is any set of primes not equal to the whole set.  We  initially consider finitely generated metabelian groups. As localization preserves surjectivity, we will see that it suffices to treat the case of  JNVN$_\pi$ groups i.e   groups which are in $\upi$ and are not virtually nilpotent but all of whose  proper quotient groups  are either  virtually nilpotent or not in $\upi$.
\begin{theorem} Let $\pi$ be a proper subset of primes and suppose that $G$ is a  finitely generated metabelian JNVN$_\pi$
group. Then $G_\pi$ is not soluble. \end{theorem}
\begin{pf}  As $G$ is a finitely generated metabelian group, $G$
satisfies max-n. Therefore we can choose $A\triangleleft G$
maximal subject to $G'\le A$ and $A$ being abelian. If $x^m \in
A,$ for some $x\in G,$ $m\in \pi',$  then $(x^m)^a=x^m $, for all
$a\in A.$  As $G\in\upi$, this implies that
$x^a=x$ and by the maximality of $A$, $x\in  A.$  Therefore $G/A$
is a $\pi'$-torsion-free group.

\par  Suppose that $A\le N\triangleleft G $ and $Z=Z(N)\not =1$,  then
$Z$ is a proper subgroup of $A$. If $x^mZ=y^mZ,$ for some $x,y \in
G,\, m\in \pi' ,$ then $a^{x^{m-1}+\dots +x+1} \in Z$ where
$a=x^{-1}y .$ As $Z\le A$ and $G/A$ is $\pi'$-torsion-free, $a\in
A$. If $g \in N ,$ then
$$ [g,a]^{x^{m-1}+\dots +x+1} = [g,a^{x^{m-1}+\dots +x+1}] =1$$
($[g^h,a^h]= [g,a^h],$ for all $h\in G$) and hence
$$(x[g,a])^m=x^m [g,a]^{x^{m-1}+\dots +x+1}=x^m.$$  This implies
that  $[g,a]=1 ,$ for all $g\in N,$  and hence $a\in Z.$
Therefore $G/Z\in\upi$ and is virtually
nilpotent. So we can find a normal subgroup $K$ containing $A$ and
of finite index in $G$ such that $K/Z$ is nilpotent of class $c,$
for some $c.$ Let $N_1=K\cap N$ and set $B_i= \gamma_i (K)\cap N_1.$ By Lemma 3.2
\cite{zh},
$$[N_1',_m K] \le \Pi_{t+s=m+2,\, t,s\ge 1}  [B_t,B_s] =1,$$ where
$m=2c-1 .$  If $Z(K)\not=1, $ then $G/Z(K)$ would be virtually
nilpotent and hence $K$ and $G$ would be virtually nilpotent
contradicting the hypothesis. Therefore $Z(K)=1$ and $N_1' =1 .$
Hence $N_1$ is abelian and, by maximality, is contained in $A$ (in fact
$A$ is  the Fitting subgroup of $G$). So that $N$ is an abelian-by-finite $\pi$-group.  If $g^n\not \in A,$ for any integer $n\ge 1,$ then let
$N=<A,g>\triangleleft G.$ By the above,  $C_A(g)=Z(N)=1.$

 \par If $G/A$ is not torsion-free we can  find $N\triangleleft G$ of finite $\pi$-index in $G$  containing $A$ such that $N/A$ is torsion-free.
  We will show that  $N_\pi$ is not soluble.  As $G_\pi$ is a  finite extension of $N_\pi$,  this would imply that  $G_\pi$
 is not soluble, by
 Theorem
 2 \cite{rod}. Therefore we can assume, without loss of generality that $G/A$ is torsion-free.

\vskip
.5cm If $A$ is not torsion-free, then, for some $p\in \pi ,$  $B=\{a\in A\mid a^p=1\}
\not= 1$.   If $x,y\in G$ and $x^mB=y^mB,$
for some $m\in \pi',$ then $y=xa,$ for some $a\in A$ with
$a^{x^{m-1} +\dots +x+1}\in B.$ Therefore
$$(xa^p)^m=x^m(a^p)^{x^{m-1} +\dots +x+1}=x^m$$ and, as $G\in \upi$, this implies that $a^p=1.$ Hence $a\in B$ and
$G/B\in\upi .$  Therefore $G$ has a normal subgroup
$K$ containing $A$ and of finite index in $G$ such that $K/B$ is
nilpotent of class $c,$ for some $c.$ If $A\not =B,$ then there
exists   a non-trivial element   $aB\in A/B\cap Z(K/B).$ Therefore
$[a,x]\in B,$ for all $x\in K. $ Hence $[a^{p}, x]=[a,x]^p=1$ and
$a^p \in Z(K)=1.$ Therefore $A$ is either torsion-free or is an
elementary abelian $p$-group, for some $p\in \pi.$
\par $\;$

   \par Suppose that  $A$ is an
elementary abelian $p$-group, for some $p\in \pi.$
If  $g\not \in
 A,$  then $C_A (g)=1.$ Suppose that $C_G (g)\cap C_G (g)^z\not=1,$
 then there exists $y,y^z\in C_G (g).$ Therefore $[x,z]\in C_G
 (y)=1$, as $y\not \in A.$ Hence $C_G (g)$ is malnormal and torsion-free,
 for all $g\not\in A.$ Therefore $G\in \ppi$ and $G_\pi$ is not soluble, by Theorem 2.6.
\par $\;$
  \par If $A$ is a torsion-free abelian group, then let $Q=G/A$ (recall we are assuming that $Q$ is torsion-free) and let $H = (A_\pi \rtimes G) /T $ where $T= <{(\po (a),
a^{-1})}\mid a\in A>^{Ncl} $. The structure of $H$ can be seen
more clearly by using the following commutative diagram:
$$\comdiag{1& \mapright{}&A&\mapright{}&G&\mapright{
}&Q&\mapright{}&1\cr&& \lmapdown{\po}&&\mapdown{}&&\Vert{}&&\cr
1&\mapright{}&A_\pi &\mapright{}& H&\mapright{ }&Q&\mapright{}&1
.}$$ As $G$ embeds in $H,$ we shall assume that $G\le H$ and that
$A\le A_\pi$. As $A_\pi$ and $Q$ are abelian, $H$ is metabelian.
  By Proposition 1.3 and Corollary 1.6  \cite{casacast}, $G_\pi \cong H_\pi
  .$

 \vskip .5cm Let $x\in H$ where $x= ga,$ for some
 $g\in G ,a\in A_\pi.$ If $ x\in H$ has no m-th root, for some
 $m\in \pi',$ then  $x\not\in A_\pi $ and $g\not\in A.$  Therefore $C_A(g)=C_A (x)=1$  and, as $A_\pi /A$ is a
$\pi'$-torsion group, this implies that
$C_{A_\pi} (x) =1$.
 If $C_H (x)\cap C_H(x)^z\not=1,$  then
  there exists $y\in C_H(x)$ such that $y^z\in C_H (x).$  Hence $ [x,z]\in C_H (y)\cap A_\pi=1,$ as $y\not\in A_\pi .$
  Therefore  $C_H (x)$ is malnormal and torsion-free.

  \vskip .5cm Suppose that $x^m=y^m,$ for some $x,y \in H,$ $m\in
  \pi'$. Since $Q$ is a torsion-free abelian group,
  $xy^{-1}\in A_\pi.$
 If $x\in A_\pi,$
 then $y\in A_\pi$ and,  as $A_\pi\in \upi $, $x=y.$
  If $x\not\in A_\pi,$
 then $xy^{-1}\in C_{A_\pi} (x^m)=1$
 and  $H\in \upi$. Therefore  $H\in \ppi$ and $G_\pi\cong H_\pi$ is not
  soluble, by Theorem 2.6 and Theorem 2 \cite{rod}. $\square$

\end{pf}
 \par We remark here that the groups in Theorem 3.1 which lie in
 $\ppi$ have elements which  satisfy condition 2.2 (ii) (a) so that 2.2 (ii) (b) may seem superfluous. However
  let
 $G\in \ppi$ be a finitely generated metabelian  JNVN$_\pi$
 group where
 $A$ is an elementary  abelian $2$-group. If $x\in G\backslash A$ has no $m$-th root, for some $m\in \pi'$,
 then adjoin $\pi'$-roots for $x$ by letting $\displaystyle K =G*_D P$ as in section 2. If   $a,\,b\in A,$
  $z^m=x$  and   $g=ab^z,$  then $g^a=b^za=g^{-1} $ and $C_K (g)\cong C_\infty,$ by Lemma 2.3.  Therefore
  $D_K (g)\cong D_\infty .$

\begin{corollary} Let $G$ be a finitely generated metabelian group and
let $\pi$ be a finite set of primes. Then
$G$ either has a normal subgroup $N$ such that $G/N\in \ppi$ and so $G_\pi$ is  not soluble or  $G_\pi$ is
virtually nilpotent.
\end{corollary}
\begin{pf} As $G$ is a finitely generated metabelian group $G$
satisfies max-n and hence either  every quotient of $G$ that is in $\upi$ is virtually nilpotent or $G$ has a JNVN$_\pi$
quotient. In the first case $G/\Ker (G\rightarrow G_\pi)$ is
virtually nilpotent and hence $G_\pi$ is virtually nilpotent, by
Theorem 3.3 \cite{casacast}. If $G$ has  a JNVN$_\pi$ quotient,
then $G_\pi$ is not soluble by Theorem 3.1. $\square$ \end{pf}
Recall that a group $G$ is a just non-virtually nilpotent group (JNVN)  if $G$ is not virtually nilpotent but all of  its proper quotients are virtually nilpotent. These groups have been studied extensively so that their structure is well-known  \cite{mdf, robwil}.
\begin{lemma} Let $G$ be a finitely generated group which is not virtually nilpotent. Then $G$ has a JNVN quotient.
\end{lemma}
\begin{pf} The proof is a minor adaptation of the proof of 7.4.1 p. 138 \cite{lenrob}. Let $\{ N_\lambda \mid \lambda\in \Lambda\}$ be a
chain of normal subgroups of $G$ such that each $G/N_\lambda$ is
not virtually nilpotent. Let  $N= \cup_{\lambda} N_\lambda .$ Then
$G/N$ is not virtually nilpotent, as otherwise  it would be
finitely presented and hence $N$ would be the normal closure in
$G$ of a finite number of elements and hence $N=N_\lambda,$ for
some $\lambda.$
  $\square$ \end{pf}

\begin{theorem} Let  $G$ be a finitely generated soluble group.
Then either $G_\pi$ is virtually nilpotent, for all sets of primes
$\pi ,$ or there exists a finite set of primes $\tau$ such that
$G_\pi$ is not soluble whenever $\tau\subseteq \pi .$
\end{theorem}
\begin{pf} Let $G$ be a finitely generated
 soluble group. Then, by Lemma 3.3, $G$
is either virtually nilpotent or has a JNVN quotient. If $G$ is
virtually nilpotent, then $G_\pi$ is also virtually nilpotent, by
Theorem 3.3 \cite{casacast}, for any set of primes $\pi .$ If $G$
has a JNVN quotient, then, as localization preserves surjectivity, we shall assume, without loss of generality, that $G$ is a JNVN group.  By Theorem 2.2 \cite{mdf}, the Fitting
subgroup, $A$, of $G$ is abelian and is either torsion-free or has
prime exponent $p.$ If $A$ is torsion-free, then, by Theorem 3.2
\cite{mdf}, $G$ is  virtually metabelian. If $A$ has prime
exponent $p,$ then $A$ is not finitely generated as otherwise $G$
would be virtually nilpotent. Therefore $G$ is a just
non-polycyclic (JNP) group and by the remarks preceding 7.4.6
\cite{lenrob}, we see that $G$ is virtually metabelian.  If $A$ is a $p$-group, then $G  $ is virtually  a
residually finite $p$-group and if $A$ is torsion-free, then
 $G$ is virtually a residually finite $p$-group, for all but a finite number of primes $p, $ by Theorem A \cite{segal1}. In either case there exists a prime $p$ and $n\in \N$ such that $H=G^n $ is a finitely generated metabelian residually finite $p$-group.  Let $\tau$ be the finite set of primes containing $p$ and the primes dividing $n.$ If  $\pi$ is any set of primes (not equal to the set of all primes) containing $\tau,$ then $H\in \upi$ (as it is residually a finite $p$-group) and is not virtually nilpotent. So, by Corollary 3.2, $H_\pi$ is not soluble. As $G_\pi$ is an  extension of $H_\pi$ by a finite $\pi$-group,  this would imply that  $G_\pi$
 is not soluble, by
 Theorem 2 \cite{rod}. $\square$
\end{pf}

\bibliography{solja}

\end{document}